\input amstex
\input epsf
\documentstyle{amsppt}

\topmatter
\def\R{{\Bbb R}}
\def\Z{{\Bbb Z}}
\def\C{{\Bbb C}}
\def\N{{\Bbb N}}

\def\GL{{\operatorname{GL}}}
\def\vf{{\frak X}}
\def\m{{\Cal M}}
\def\mot{{\m^{12}}}

\def\fl{{\Cal L}}
\def\fly{{\fl_{\gamma}}}
\def\flo{{\fl_o}}

\def\cfly{{\widetilde\fly}}
\def\cflo{{\widetilde\flo}}
\def\p{{\Cal P}}

\def\py{\p_{\gamma}}

\def\cpy{\tilde\py}

\def\pyo{\py^1}
\def\pyt{\py^2}
\def\pyeven{{\py^{\text{even}}}}
\def\pyodd{{\py^{\text{odd}}}}
\def\cpyo{\cpy^1}
\def\cpyt{\cpy^2}
\def\ddt{{\tfrac\partial{\partial t}}}
\def\dds{{\tfrac\partial{\partial s}}}
\def\ddsu{{\tfrac{\partial u}{\partial s}}}
\def\ddtu{{\tfrac{\partial u}{\partial t}}}
\def\img{{\operatorname{img}}}

\def\Hom{{\operatorname{Hom}}}
\def\ev{{\operatorname{ev}}}
\def\resp{{resp.\ }}
\def\ie{{i.e.\ }}
\def\cf{{cf.\ }}

\def\vfbyreg{\vf^\beta_{\text{$\gamma$-reg}}}
\def\vfboreg{\vf^\beta_{\text{$\gamma_0$-reg}}}
\def\Jreg{{\Cal J}_{\text{reg}}}
\def\id{\operatorname{id}}
\def\Wh{\operatorname{Wh}}
\def\lt{\operatorname{lt}}

\def\bbs{{\bold b}^*}
\def\bcs{{\bold c}^*}
\def\bds{{\bold d}^*}
\def\bhs{{\bold h}^*}
\def\bbss{{\bold b}^{*+1}}
\def\bcss{{\bold c}^{*+1}}
\def\bcsss{{\bold c}^{*+2}}
\def\grad{\operatorname{grad}}

\title
Non-contractible periodic trajectories of symplectic vector fields,
Floer cohomology and symplectic torsion.
\endtitle
\rightheadtext{Floer cohomology and non-contractible periodic orbits \dots}

\dedicatory
(preliminary version)
\enddedicatory

\author
Dan Burghelea and Stefan Haller
\endauthor

\address
The Ohio State University, 231 W 18th Avenue, Columbus, OH 43210 USA
\endaddress

\email
burghele\@mps.ohio-state.edu
\endemail

\address
University of Vienna, Strudlhofgasse 4, A-1090 Vienna, Austria
\endaddress

\email
stefan\@mat.univie.ac.at
\endemail

\thanks
The first author was supported by NSF.
\endthanks

\thanks
The second author is supported by the Fonds zur F\"orderung der 
wissenschaftlichen Forschung' (Austrian Science Fund), project number P14195
\endthanks

\abstract
For a closed symplectic manifold $(M,\omega)$, a compatible
almost complex structure $J$, a $1$-periodic time dependent symplectic
vector field $Z$ and a homotopy class of closed curves $\gamma$
we define a Floer complex based on $1$-periodic trajectories of $Z$ in the
homotopy class $\gamma$. We suppose that the closed $1$-form $i_{Z_t}\omega$
represents a cohomology class $\beta(Z):=\beta$, independent of $t$.
We show how to associate to $(M,\omega,\gamma,\beta)$ and to two pairs 
$(Z_i,J_i)$, $i=1,2$ with $\beta(Z_i)=\beta$ an invariant, the
relative symplectic torsion, which is an element in the Whitehead group
$\Wh(\Lambda_0)$, of a Novikov ring $\Lambda_0$ associated with
$(M,\omega,Z,\gamma)$. If the cohomology of the Floer complex vanishes or if 
$\gamma$ is trivial we derive an invariant, the symplectic torsion for any 
pair $(Z,J)$.

We prove, that when $\beta(\gamma)\neq 0$, or when $\gamma$ is non-trivial
and $\beta$ is \lq small\rq,
the cohomology of the Floer complex is trivial, but the symplectic torsion 
can be non-trivial. Using the first fact we conclude results about
non-contractible $1$-periodic trajectories of $1$-periodic symplectic
vector fields.
In this version we will only prove the statements for closed weakly monotone
manifolds, but note that they remain true as formulated for arbitrary
closed symplectic manifolds.
\endabstract

\toc
\widestnumber\head{9}
\head 1. Introduction\endhead
\head 2. Topological data and associated Novikov rings\endhead
\head 3. Analytical data\endhead
\head 4. The Floer complex\endhead
\head 5. Floer cohomology and the proof of Theorem 1.1\endhead
\head 6. Symplectic torsion\endhead
\head 7. Proof of Proposition 6.5\endhead
\head 8. Appendix (An example)\endhead
\endtoc
\endtopmatter

\document

\head
1. Introduction
\endhead

Let $(M,\omega)$ be a closed connected symplectic manifold, and
$Z\in C^\infty(S^1,\vf(M,\omega))$ a $1$-periodic symplectic vector field on
$M$. We are interested in $1$-periodic trajectories of $Z$. If every $Z_t$
is Hamiltonian, Floer theory (initially developed by Floer,
Hofer-Salamon for weakly monotone symplectic manifolds and later extended by
Liu-Tian, Fukaya-Ono to all symplectic manifolds) was used by Salamon-Zehnder
\cite {SZ92}
to count the $1$-periodic trajectories that are contractible.
In fact, when all these trajectories are non-degenerate
they generate a cochain complex which computes $H^*(M;\Z)$. As a consequence
the number of non-degenerate $1$-periodic contractible trajectories can be
estimated from below by the sum of the Betti numbers of $M$. In \cite{LO95}
this was extended to arbitrary $1$-periodic time dependent symplectic
vector field $Z$ as above.

The discussion in the above mentioned papers is done under the additional
hypotheses, that the symplectic manifold $(M,\omega)$ is weakly monotone, see
Definition~4.3, but in view of the work of Liu-Tian \cite{LT98} and
Fukaya-Ono \cite{FO99} this hypothesis can be removed.

We will go one step further and consider $1$-periodic trajectories which are
non-contractible under the hypothesis that 
$\beta:=[i_{Z_t}\omega]\in H^1(M;\R)$
does not depend on $t$. We consider:
\roster
\item
a homotopy class $\gamma$ of closed curves,
\item
a pair $(Z,J)$, where $Z$ is a $1$-periodic symplectic vector field whose
$1$-periodic trajectories in $\gamma$ are all non-degenerate and such
that $\beta:=[i_{Z_t}\omega]\in H^1(M;\R)$ is independent of $t$, $J$ is a
compatible almost complex structure (\cf section~4) and
\item
additional data, consisting of a c-structure above $\gamma$
(\cf Definition~2.1) and `coherent orientations'.
\endroster
In this paper we suppose for simplicity that $(M,\omega)$ is weakly monotone
and the pair $(Z,J)$
satisfies a regularity condition with respect to $\gamma$ which is a generic property
(\cf section~4). In view of the work \cite {LT98} and \cite{FO99} these two hypotheses
are not necessary if the Novikov ring associate to
$(M,\omega, \beta, \gamma)$ has coefficients containing the rational numbers
\footnote{ in view of more recent work \cite {FO00} not necessary at all.}.
Under the hypotheses of weak monotonicity and $\gamma$-regularity the concept of
`coherent orientations' is the one considered in \cite{FH93}. Without these hypotheses
`coherent orientations' should be understood as `compatible Kuraninishi structures with
corners' as discussed in \cite{FO99} section~19.

In sections~3 and 4 we associate to this data a Floer complex, which is a
cochain
complex of $\Z_{2N}$-graded free $\Lambda$-modules. Here $\Lambda$ is a
Novikov ring associated to $(M,\omega,\beta,\gamma)$, \cf section~2, and the
integer $N$ is the minimal Chern number as defined in section~2, depending on
$(M,\omega,\gamma)$.

The underlying free module of this complex is generated by the
$1$-periodic closed trajectories of $Z$ in the class $\gamma$, which are
non-degenerate in view of the hypotheses on $Z$,
hence are finitely many.

As expected one shows that the cohomology of this complex is independent of
the pair $(Z,J)$. If $\beta(\gamma)\neq 0$, or if $\gamma$ is non-trivial and
$\beta$ is \lq small\rq, this cohomology vanishes.
These properties are collected in Theorems~5.7 and 5.8 of which
one derives the following geometric result. It remains true
as stated without the weak monotonicity hypothesis.

\proclaim{Theorem 1.1}
Let $(M,\omega)$ be a closed connected weakly monotone symplectic
manifold and $Z$ a $1$-periodic time dependent symplectic vector field,
such that the cohomology class $\beta=[i_{Z_t}\omega]\in H^1(M;\R)$
is independent of $t$. Then:
\roster
\item
If in a given non-trivial homotopy class of closed curves $\gamma$, so that 
$\beta(\gamma)\neq 0$, there exists a
non-degenerate $1$-periodic trajectory of $Z$, then there exits at least
one more geometrically different, possibly degenerate trajectory in the same
homotopy class.
\item
If all $1$-periodic trajectories of $Z$ in a given non-trivial homotopy
class $\gamma$, so that $\beta(\gamma)\neq 0$, are non-degenerate, then their 
number is even. 
\endroster
\endproclaim

\remark{Remark 1.2}
This theorem remains true without the hypothesis $\beta(\gamma)\neq 0$.
K.~Ono has informed us, that he can also prove this result.
\endremark


Our next observation is, that the cochain complex associated to
$(M,\omega,\gamma,Z,J)$ and the additional data has a well defined torsion,
independent of the additional data, at least in the case that 
$\beta(\gamma)\neq 0$ or $\gamma$ trivial. This torsion, a priori defined in
$\Wh(\Lambda)$, actually lies in the subgroup $\Wh(\Lambda_0)$. The
rings $\Lambda$ and $\Lambda_0$ are defined in section~2 and the groups
$\Wh(\Lambda)$ and $\Wh(\Lambda_0)$ in section~6.

We hope that this symplectic torsion will be an important new invariant
to be used in the study of $1$-periodic closed trajectories of $Z$ and maybe
other Floer complex related problems, particularly when the Floer complex is
acyclic. At present we know little about it (despite of a good number of
conjectures we have). For example we know that
it is non-trivial, \cf section~8, but we do not know yet how it depends on
the almost complex structure $J$. We also know how to calculate this
torsion in the case $\gamma$ is trivial and $Z$ is time independent.
More precisely, consider the Riemannian metric
$g$ induced from $\omega$ and $J$ \cf section~4. The symplectic torsion can be expressed
both in terms of the dynamics of the vector field $\grad_g(i_Z\omega)$, precisely its closed
trajectories, and in terms of the spectral theory of the (time dependent) Laplacians
of the complex $(\Omega^*(M),d^*+t(i_Z\omega)\wedge\cdot)$ with respect
to the metric $g$.  The second expression requires some
Dirichlet series theory as in \cite{BH01}.
When the cohomology of $M$ with coefficients in $\beta=[i_Z\omega]$ is
trivial our symplectic torsion identifies to
the `zeta function', defined in \cite{Hu00}.

We point out that section~1.9 of \cite{Hu00}, entitled
`Possible generalizations' refers to ideas close to those presented here.
We thank A.~Fel'shtyn for bringing this paper to our attention.
We will return to this symplectic torsion in future research.
We thank K.~Ono for pointing out a mistake in a previous version of this
paper. We understand, that he was also aware of the possibility to define a
Floer complex for non-contractible trajectories.


\head
2. Topological data and associated Novikov rings
\endhead

Let $(M,\omega)$ be a closed connected symplectic manifold of dimension $2n$
and let $\fl:=C^\infty(S^1,M)$ denote the space of
smooth free loops, where we think of $S^1$ as $\R/\Z$. The first Chern
class $c_1\in H^2(M;\Z)$ of $(M,\omega)$ defines a cohomology class
$\bar c_1\in H^1(\fl;\Z)=\Hom\big(H_1(\fl;\Z),\Z\big)$, by interpreting a
class in $H_1(\fl;\Z)$ as linear combination of `tori', \ie maps
$S^1\times S^1\to M$, and integrating $c_1$ over it. Similarly
$[\omega]\in H^2(M;\R)$ gives rise to a cohomology
class $\overline{[\omega]}\in H^1(\fl;\R)=\Hom\big(H_1(M;\Z),\R\big)$. Given a
cohomology class $\beta\in H^1(M;\R)$ we define
$$
\phi:=\overline{[\omega]}+\ev^*\beta\in H^1(\fl;\R),
$$
where $\ev:\fl\to M$ denotes the evaluation of a loop at the point
$0\in S^1=\R/\Z$.

\definition{Definition~2.1}
A c-structure for $(M,\omega)$ is a connected component of the space of
$\R$-vector bundle homomorphisms $\tilde y:S^1\times\C^n\to TM$, satisfying
$\omega\big(\tilde y_t(v),\tilde y_t(iv)\big)>0$, for all $t\in S^1$ and all
$0\neq v\in\C^n$.
\footnote{Notice that every such homomorphism is a fiber wise isomorphism.}
We denote by $[\tilde y]$ the c-structure represented by $\tilde y$.
\enddefinition

By assigning to $\tilde y$ the underlying map $y:S^1\to M$ one obtains a
surjective mapping $\pi:\{\text{c-structures}\}\to\pi_0(\fl)$,
and we say that $[\tilde y]$ is a c-structure above
$\gamma=\pi([\tilde y])$. The group $\pi_1(\GL_n(\C))\cong\Z$ acts freely on
the set of c-structures with $\pi_0(\fl)$ as the set of orbits.
By fixing a c-structure $[\tilde y]$ we specify a homotopy class of closed
curves $\gamma$ or equivalently a component $\fl_\gamma$ of $\fl$. By choosing
a representative $\tilde y$ of $[\tilde y]$ we specify a base point $y$ in
$\fl_\gamma$.

Note, that if $\gamma$ is trivial, \ie $\fly$ is the component of contractible
loops, then there exists a canonic c-structure above $\gamma$, given by any
$\tilde y:S^1\times\C^n\to TM$, which is constant in $t\in S^1$.

We choose a c-structure $[\tilde y]$ and set $\gamma=\pi([\tilde y])$.
Consider the connected Abelian principal covering $\pi:\cfly\to\fly$, with
structure group
$$
\Gamma:=\frac{\pi_1(\fl,y)}{\ker\bar c_1\cap\ker\phi},
\tag{2.1}
$$
where $\bar c_1$ and $\phi$ are considered as homomorphisms
$\pi_1(\fl,y)\to\R$. Clearly
$\pi^*\phi=\pi^*\bar c_1=0\in H^1(\cfly;\R)$, and $\bar c_1$ and $\phi$
induce homomorphisms
$$
\bar c_1:\Gamma\to\Z\quad\text{and}\quad\phi:\Gamma\to\R.
$$

We define the minimal Chern number $N\in\N$ by
$$
N\Z=\img\big(H_1(\fl_\gamma;\Z)@>{\bar c_1}>>\Z\big)\subseteq\Z,
$$
with the convention $N=\infty$ if this image is $0$. $N$ depends only on
$(M,\omega,\gamma)$. If $\gamma$ is trivial this is the usual minimal Chern
number as discussed in \cite{HS95}, but in general it is smaller.
Indeed, the diagram
$$
\CD
H_1(\fl_\gamma;\Z) & @>{\bar c_1}>> & \Z
\\
@AAA & & @AA{c_1}A
\\
\pi_1(\fl,y) & @<<< & \pi_2(M)
\endCD
$$
commutes, where the bottom mapping
$\pi_2(M)=\pi_1(\Omega(M,y(0))\to\pi_1(\fl,y)$ is induced from the mapping
$\Omega(M,y(0))\to\fl$ given by concatenating a loop based at $y(0)$ with $y$,
and $\Omega(M,y(0))$ denotes the space of base pointed loops.

Finally we choose a commutative ring with unit $R$ and let
$\Lambda=\Cal N(\Gamma,\phi,R)$ denote the Novikov ring, associated to
$\Gamma$ and the weighting homomorphism $\phi:\Gamma\to\R$ with values in
$R$. More precisely $\Lambda$ consists of all functions $\lambda:\Gamma\to R$,
such that
$$
\big\{A\in\Gamma\bigm| \phi(A)\leq c,\lambda(A)\neq 0\big\}
$$
is finite for all $c\in\R$, with multiplication
$$
(\lambda*\nu)(A)=\sum_{B\in\Gamma}\lambda(B)\nu(B^{-1}A).
$$
$\Lambda$ is a commutative ring with unit, which depends on
$(M,\omega,\beta,\gamma,R)$. We also consider the group
$$
\Gamma_0:=\ker\bigl(\bar c_1:\Gamma\to\Z\bigr)
$$
and the Novikov ring
$\Lambda_0=\Cal N(\Gamma_0,\phi,R)$ associated to $\phi:\Gamma_0\to\R$.
Clearly $\Lambda_0$ is a subring of $\Lambda$. If $R$ has no zero divisors
then $\Lambda$ and $\Lambda_0$ are rings without zero divisors, for
$\Gamma$ and $\Gamma_0$ are free Abelian. If $R$ is a principal
ideal domain \resp a field, then so is $\Lambda_0$, since
$\phi:\Gamma_0\to\R$ is injective, \cf~\cite{HS95}

\head
3. Analytical data
\endhead

Recall that $\fl$ is a Fr\'echet manifold with tangent bundle
$C^\infty(S^1,TM)$. So a tangent vector at $x\in\fl$ is simply a vector field
along $x$. For $\sigma\in\Omega^k(M)$ we define $\tilde\sigma\in\Omega^k(\fl)$
by
$$
\tilde\sigma(X^1_x,\dotsc,X^k_x):=
\int_{S^1}\sigma\big(X^1_x(t),\dotsc,X^k_x(t)\big)dt,
$$
where $X^i_x$ are vector fields along the common loop $x$. Every
$\tilde\sigma$ is $S^1$-invariant, \ie $L_\zeta\tilde\sigma=0$, where
$\zeta$ is the fundamental vector field
\footnote{
Recall that if $\mu: S^1\times X\to X$ is a smooth action, then the fundamental
vector field $\zeta\in\vf(X)$ is given by $\zeta(x):=\ddt|_{t=0}\mu(t, x)$.
}
of the natural $S^1$-action on $\fl$, and one has 
$d\tilde\sigma=\widetilde{d\sigma}$. Therefore
$$
\Omega^*(M)\to\Omega^{*-1}(\fl),\quad
\sigma\mapsto\bar\sigma:=(-1)^{|\sigma|+1}i_\zeta\tilde\sigma
$$
commutes with the differentials and every $\bar\sigma$ is $S^1$-invariant, too.
Moreover, for $f:N\to\fl$ and $\sigma\in\Omega^*(M)$ one has
$$
\int_Nf^*\bar\sigma=\int_{N\times S^1}\hat f^*\sigma,
\tag{3.1}
$$
where $\hat f:N\times S^1\to M$, $\hat f(z,t)=f(z)(t)$. Because of \thetag{3.1} the
induced mapping $H^*(M;\R)\to H^{*-1}(\fl;\R)$ maps the first Chern class of
$(M,\omega)$ to the class $\bar c_1$ we have defined in section~2. Moreover
$\tilde\omega$ is an $S^1$-invariant closed weakly non-degenerate $2$-from on
$\fl$. From \thetag{3.1} we also see, that $\bar\omega$ is a form representing
$\overline{[\omega]}\in H^1(\fl;\R)$ from section~2.

Let $Z\in C^\infty(S^1,\vf(M))$ be a 1-periodic vector field. It defines a
vector field $\tilde Z$ on $\fl$, by setting $\tilde Z_x(t):=Z_t(x(t))$, for
$x\in\fl$ and $t\in S^1$. It is straightforward to show

\proclaim{Lemma 3.1}
If $Z$ is symplectic, \ie $di_{Z_t}\omega=L_{Z_t}\omega=0$ for all $t\in
S^1$,
then $\tilde Z$ is symplectic, \ie $di_{\tilde Z}\tilde\omega=L_{\tilde
Z}\tilde\omega=0$. If moreover $[i_{Z_t}\omega]=\beta\in H^1(M;\R)$ for all
$t\in S^1$ then $[i_{\tilde Z}\tilde\omega]=\ev^*\beta\in H^1(\fl;\R)$.
\endproclaim

Suppose $Z$ is symplectic and define the action 1-form
$$
\alpha:=i_{\tilde Z}\tilde\omega+\bar\omega=i_{\tilde Z-\zeta}\tilde\omega
\in\Omega^1(\fl).
\tag{3.2}
$$
By Lemma~3.1 $\alpha$ is closed, \ie $\tilde Z-\zeta$ is a symplectic
vector field
on $\fl$. The following is an immediate consequence of the weak non-degeneracy
of $\tilde\omega$.

\proclaim{Lemma 3.2}
The zeros of $\alpha$ are precisely the 1-periodic solutions $x:S^1\to M$ of
$$
x'(t)=Z_t(x(t)).
\tag{3.3}
$$
\endproclaim

If $[i_{Z_t}\omega]=\beta\in H^1(M;\R)$ for all $t\in S^1$, then the second
part of Lemma~3.1 says, that $\alpha$ is a closed one form representing
$\phi\in H^1(\fl;\R)$. From \thetag{2.1} we then see, that
$\pi^*\alpha=da$ for some $a\in C^\infty(\cfly,\R)$ and
$$
a(A\sharp\tilde x)=a(\tilde x)+\phi(A),
$$
for $\tilde x\in\cfly$ and $A\in\Gamma$.
\footnote{Note that $a$ is unique up to a additive constant.}

For $\gamma\in\pi_0(\fl)$ let $\py$ denote the set of $1$-periodic trajectories
of \thetag{3.3} in $\fly$ and set $\cpy:=\pi^{-1}(\py)$, where
$\pi:\cfly\to\fly$ is the covering. If all $x\in\py$ are
non-degenerate, that is the corresponding fixed points $x(0)\in M$ of the
time $1$ flow $\Psi^Z_1$ to $Z$ are non-degenerate, then $\py$ is finite.

We denote by $\vf^\beta$ the $1$-periodic time dependent symplectic
vector fields $Z\in C^\infty\bigl(S^1,\vf(M)\bigr)$, such that
$[i_{Z_t}\omega]=\beta\in H^1(M;\R)$, for all $t\in S^1$. Moreover we denote
by $\vfbyreg$ the subset of $\vf^\beta$ consisting of vector fields whose
$1$-periodic trajectories in $\fly$ are all non-degenerate.

%
%

\proclaim{Proposition 3.3}
Let $\beta$ and $\gamma$ be as above. Then the set $\vfbyreg$ is a generic
subset of $\vf^\beta$, in the sense of Baire.
\footnote{This is meant in the same sense, as in \cite{HS95}.}
\endproclaim

The genericity of $\vfbyreg$ is essentially  proven in Theorem~3.1 in 
\cite{HS95} and Theorem~3.1 in \cite{LO95}.


Let $[\tilde y]$ be a c-structure set $\gamma=\pi([\tilde y])$ and suppose
$Z\in\vfbyreg$. The c-structure $[\tilde y]$ gives rise to a
well defined, up to homotopy, $\C$-vector bundle trivialization of $x^*TM$ for every
$x\in\fly$. So one gets a well defined Conley-Zehnder index
$\mu^{[\tilde y]}:\cpy\to\Z$, satisfying
$$
\mu^{[\tilde y]}(A\sharp\tilde x)=\mu^{[\tilde y]}(\tilde x)+2\bar c_1(A),
$$
for all $\tilde x\in\cpy$ and $A\in\Gamma$, \cf \cite{SZ92}.

\remark{Remark 3.4}
If one changes the c-structure $[\tilde y]$ by an element in
$\pi_1(\GL_n(\C))\cong\Z$ (\cf section~2) the component $\fly$ remains the
same and $\mu^{[\tilde y]}$ experiences a shift by a constant in $2\N$.
In particular for $\tilde x_-,\tilde x_+\in\cpy$ the difference
$\mu^{[\tilde y]}(\tilde x_+)-\mu^{[\tilde y]}(\tilde x_-)$
depends only on $\gamma$. We will drop $[\tilde y]$ from the notation and
simply write $\mu$ for the index.
We also have induced index maps $\mu:\cpy\to\Z_{2N}$ and $\mu:\cpy\to\Z_2$.
Obviously the last one does not depend on the c-structure over $\gamma$.
\endremark

Let $C_F^*$ denote the set of all functions $\xi:\cpy\to R$, such that
$$
\big\{\tilde x\in\cpy\bigm| a(\tilde x)\leq c,\xi(\tilde x)\neq 0\big\}
$$
is finite for all $c\in\R$. $C_F^*$ is a free
$\Z_{2N}$-graded $\Lambda$-module, via
$$
(\lambda*\xi)(\tilde x)=\sum_{A\in\Gamma}\lambda(A)\xi(A^{-1}\sharp\tilde x).
$$
The component $C^i_F$, $i\in\Z_{2N}$, consists of the functions $\xi$ which
vanish on $\tilde x$ with $\mu(\tilde x)\neq i$. The total rank of $C^*_F$
equals the number of $1$-periodic trajectories of $Z$ in $\fly$, \ie the
cardinality of $\py$. $C^*_F$ depends on $(M,\omega,[\tilde y],Z,R)$, but
the associated $\Z_2$-graded module does only depend on
$(M,\omega,\gamma,Z,R)$, for a change of the c-structure over $\gamma$
shifts the grading of $C^*_F$ by an even integer.

\head
4. The Floer complex
\endhead

A smooth almost complex structure $J$ on $M$ is called compatible with
$\omega$ if
$$
g(X,Y):=\omega(X,JY)
$$
defines a Riemannian metric on $M$.
We will denote by $\Cal J$ the set of
almost complex structures compatible with $\omega$.

For a closed symplectic manifold $(M,\omega)$, a c-structure $[\tilde y]$
above $\gamma$ and a pair $(Z,J)\in\vfbyreg\times\Cal J$ one can associate
a cochain complex of free $\Lambda=\Cal N(\Gamma,\phi,R)$ modules. For reason
of simplicity we will do this only under the additional hypotheses that
$(M,\omega)$ is weakly monotone and the pair is $\gamma$-regular, a
concept defined below.
This is not an inconvenient restriction for the study of the closed
trajectories of $Z$ because in view of Proposition~4.5\therosteritem2 for
any $Z\in \vfbyreg$ there exists a $\gamma$-regular pair $(Z',J)$
with $Z$ and $Z'$ having the same trajectories.

Choose $J\in\Cal J$ and let $\tilde g$ denote the induced weak Riemannian
metric on $\fl$, given by
$$
\tilde g(X_x,Y_x)=\int_{S^1}g\bigl(X_x(t),Y_x(t)\bigr)dt.
$$
For $Z\in\vf^\beta$ we have the action $1$-form $\alpha$ \thetag{3.2},
whose gradient is a well defined vector field on $\fl$ equal to $\tilde
J\tilde Z$, where $\tilde J$ is the almost complex structure on $\fl$, $(\tilde
JX_x)(t):=JX_x(t)$. A gradient flow line is a mapping $u:\R\times S^1\to M$
satisfying
$$
D_{Z,J}u:=\ddsu+J\ddtu-JZ_t(u(s,t))=0.
\tag{4.1}
$$
For every curve $u:\R\to\fl$ we have the energy
$$
E(u):=\int_\R|u'(s)|^2_{\tilde g}ds
=\int_{\R\times S^1}|\ddsu|^2_gds dt.
$$
Using the non-degeneracy of the 1-periodic trajectories in $\py$ one shows,
\cf \cite{F89} and \cite{SZ92}:

\proclaim{Proposition 4.1}
Let $Z\in\vfbyreg$, $J\in\Cal J$ and suppose $u:\R\to\fly$ satisfies
\thetag{4.1}. Then the following are equivalent:
\roster
\item
$u:\R\to\fly$ has finite energy.
\item
There exist $x_-,x_+\in\py$, such that
$$
\lim_{s\to\pm\infty}u(s,t)=x_\pm(t)
\quad\text{and}\quad
\lim_{s\to\pm\infty}\dds u(s,t)=0
$$
both uniformly in $t\in S^1$ and exponentially in $s\in\R$.
\endroster
\endproclaim
For $\tilde x_-,\tilde x_+\in\cpy$ let $\m(\tilde x_-,\tilde x_+)$ denote the
space of finite energy solutions of \thetag{4.1} connecting $\tilde x_-$ with
$\tilde x_+$. If we want to emphasize the dependence on $(Z,J)$ we will
write $\m(\tilde x_-,\tilde x_+,Z,J)$.
For $u\in\m(\tilde x_-,\tilde x_+)$ one has
$$
E(u)=a(\tilde x_+)-a(\tilde x_-).
\tag{4.2}
$$

\definition{Definition 4.2}
A compatible almost complex structure $J\in \Cal J$ is called regular if for any
simple $J$-holomorphic curve $v:S^2\to M$ the linearization of the
Cauchy-Riemann operator $\bar\partial_J$ at $v$ is surjective, \cf
section~2 in \cite{HS95} for definitions.
\enddefinition

We will denote by $\Jreg\subseteq\Cal J$ the set of all
regular almost complex structures. It is shown in \cite{HS95} Theorem~2.2,
that $\Jreg$ is a generic subset of $\Cal J$ in the sense of Baire.

\definition{Definition 4.3}
A symplectic manifold is called weakly monotone \cite{HS95},
if for all $\tau\in\pi_2(M)$
$$
3-n\leq c_1(\tau)<0\quad\Rightarrow\quad\omega(\tau)\leq 0,
$$
where $\dim(M)=2n$.
\enddefinition
The relevance of this concept comes from the fact, that
a weakly monotone manifold, when equipped with $J\in\Jreg$, has
no $J$-holomorphic spheres of negative Chern number.
Note that every symplectic manifold of dimension smaller or equal to $6$ is
weakly monotone.

It is shown in \cite{HS95}, Proposition~2.4, that for a closed weakly monotone
symplectic manifold and $J\in\Jreg$ the subset $M_0(\infty,J)$ \resp
$M_1(\infty,J)$, consisting of points of $M$ which lie on a
non-constant $J$-holomorphic curve $v:S^2\to M$ with $c_1(v)\leq 0$ \resp
$c_1(v)\leq 1$ has  codimension greater or equal to $4$
\resp $2$, so for generic pairs $(Z,J)$, both the closed trajectories in the
homotopy class $\gamma$ and the image in $M$ of the relevant maps
$u\in\Cal M(\tilde x_-,\tilde x_+)$ have empty intersections with the sets
$M_0(\infty,J)$ and $M_1(\infty,J)$. This is a geometrically convenient
feature which permits us to define the Floer complex as in \cite{HS95},
but in view of the recent work \cite{LT98} and \cite{FO99}, not necessary.

\definition{Definition 4.4}
Let $(M,\omega)$ be a closed weakly monotone symplectic manifold and
$\gamma\in\pi_0(\fl)$. A pair $(Z,J)\in\vfbyreg\times\Jreg$ is called
$\gamma$-regular if the following conditions hold:
\roster
\item
For any $\tilde x_-,\tilde x_+\in\cpy$ and any $u\in\m(\tilde x_-,\tilde x_+)$
the linearization of the Cauchy-Riemann operator $D_{Z,J}$ at $u$ is
surjective, \cf \cite{HS95} for definitions.
\item
For any $x\in\py$ and any $t\in S^1$ one has $x(t)\notin M_1(\infty,J)$.
\item
For any $\tilde x_-,\tilde x_+\in\cpy$ with
$\mu(\tilde x_+)-\mu(\tilde x_-)\leq 2$, any $u\in\m(\tilde x_-,\tilde x_+)$
and any $(s,t)\in\R\times S^1$ one has $u(s,t)\notin M_0(\infty,J)$.
\endroster
\enddefinition

As in \cite{F89}, see also the proof of Theorem~3.2 in \cite{HS95},
\cite{Mc90}, \cite{SZ92} and \cite{LO95} one shows

\proclaim{Proposition 4.5}
Let $(M,\omega)$ be a closed weakly monotone symplectic manifold.
\roster
\item
The set of $\gamma$-regular pairs is a generic subset of
$\vf^\beta\times\Cal J$, in the sense of Baire.
\item
For $Z\in\vfbyreg$ one can find a $\gamma$-regular pair $(Z',J)$, such that
the $1$-periodic trajectories of $Z$ and $Z'$ in $\fly$ are the same, and the
vector fields $Z$ and $Z'$ agree on these trajectories up to order $2$.
Especially the Conley-Zehnder indices agree, too.
\item
For any $\gamma$-regular pair $(Z,J)$ and any $\tilde x_-,\tilde x_+\in\cpy$
the space $\m(\tilde x_-,\tilde x_+)$ is an orientable smooth manifold of
dimension $\mu(\tilde x_+)-\mu(\tilde x_-)$.
\footnote{
Although $\mu$ depends on the c-structure $[\tilde y]$ over $\gamma$, this
difference depends only on $\gamma$, \cf Remark~3.4.}
It admits a natural $\R$-action given by reparameterization, which is free
and proper if the index difference is bigger than $0$.
\endroster
\endproclaim

Suppose $(M,\omega)$ is weakly monotone and suppose $(Z,J)$ is a
$\gamma$-regular pair. The uniform bound on the energy \thetag{4.2} together
with the fact, that due to weak monotonicity no bubbling can occur,
yield (\cf \cite{F88}, \cite{S90}, \cite{HS95}, \cite{O95} and \cite{LO95}):

\proclaim{Proposition 4.6}
Suppose $\tilde x_-,\tilde x_+\in\cpy$ and $\mu(\tilde x_+)-\mu(\tilde x_-)=1$.
Then the set $\m(\tilde x_-,\tilde x_+)/\R$ is a finite. Moreover
$$
\bigl\{A\in\Gamma\bigm|\phi(A)\leq c,c_1(A)=0,\m(\tilde x_-,A\sharp\tilde
x_+)\neq\emptyset\bigr\}
$$
is finite for all $c\in\R$.
\endproclaim

As in \cite{FH93} we choose coherent orientations $\Cal O$ of
$\m(\tilde x_-,\tilde x_+)/\R$, and for $\mu(\tilde x_+)-\mu(\tilde x_-)=1$
we set
$$
n(\tilde x_-,\tilde x_+)
:=\sharp\bigl(\m(\tilde x_-,\tilde x_+)/\R\bigr)\in\Z,
$$
where the points are counted with signs according to their orientation.

If the index difference is different from $1$ we set 
$n(\tilde x_-,\tilde x_+)=0$. The weak monotonicity also implies that for
$\mu(\tilde x_+)-\mu(\tilde x_-)=2$ the one dimensional manifolds
$\m(\tilde x_-,\tilde x_+)/\R$ are compact up to broken trajectories.
Together with a gluing argument the coherent orientations yield
$$
\sum_{\tilde x\in\cpy}n(\tilde x_-,\tilde x)n(\tilde x,\tilde x_+)=0
\tag{4.3}
$$
and
$$
n(A\sharp\tilde x_-,A\sharp\tilde x_+)=n(\tilde x_-,\tilde x_+),
\tag{4.4}
$$
for all $\tilde x_-,\tilde x_+\in\cpy$ and all $A\in\Gamma$.

From Proposition~4.6 we see, that
$$
\partial(\delta_{\tilde x})
:=\sum_{\tilde x_+\in\cpy}n(\tilde x,\tilde x_+)\delta_{\tilde x_+}\in
C_F^*,
$$
and \thetag{4.4} gives $\partial(\delta_A*\delta_{\tilde
x})=\delta_A*(\partial\delta_{\tilde x})$. Here $\delta_{\tilde x}\in C_F^*$
denotes the element being 1 at $\tilde x$ and zero elsewhere. So $\partial$
extends uniquely to a $\Lambda$-linear map
$$
\partial:C_F^*\to C_F^{*+1}
$$
of degree $1$, and \thetag{4.3} immediately gives $\partial^2=0$.

To a closed weakly monotone symplectic manifold $(M,\omega)$,
a c-structure $[\tilde y]$ over $\gamma$, a symplectic vector field
$Z\in\vfbyreg$, a compatible almost complex structure $J\in\Jreg$, such that
the pair $(Z,J)$ is $\gamma$-regular, coherent orientations $\Cal O$ and a
commutative ring of coefficients $R$ we have associated a $\Z_{2N}$-graded
Floer cochain complex $C^*_F$ of free $\Lambda$-modules.
We denote the corresponding Floer cohomology $\Lambda$-module by $H_F^*$.
A priori it depends on $(M,\omega,[\tilde y],Z,J,\Cal O,R)$.

\remark{Remark 4.7}
As already observed, the Floer cochain complex can be defined even when
$(M,\omega)$ is not weakly monotone and the pair
$(Z,J)\in\vfbyreg\times\Cal J$ not necessary $\gamma$-regular. In this case,
if the pair $(Z,J)$ does not satisfies Definition~4.4\therosteritem1
(which insures that $\Cal M(\tilde x_-,\tilde x_+)$ are smooth orientable
manifolds), `coherent orientations' on $\Cal M(\tilde x_-,\tilde x_+)$ have
to be replaced by `compatible oriented Kuranishi structures with corners' on
$\tilde\Cal M(\tilde x_-,\tilde x_+)$ in the sense of Fukaya-Ono \cite{FO99}.
Here $\tilde\Cal M(\tilde x_-,\tilde x_+)$ denotes the space of stable
trajectories which is compact, Hausdorff and has Kuranishi structures,
\cf \cite{FO99}, section~19.
\endremark

\head
5. Floer cohomology and the proof of Theorem~1.1
\endhead

In this section we want to check that the cohomology of the Floer complex
depends only on $(M,\omega,\beta,[\tilde y],R)$. We will sketch the proof
only in the weakly monotone case, but the statement is true in general,
\cf Remark~5.11 below.

Suppose we have two  pairs $(Z^i,J^i)\in\vfbyreg\times\Cal J$. A homotopy between
these pairs is a pair $(Z,J)$, where $Z\in C^\infty(\R,\vf^\beta)$ and
$J\in C^\infty(\R,\Cal J)$, such that
$$
\align
(Z_s,J_s)=(Z^1,J^1)
\quad&\text{for all $s\leq -1$ and}
\\
(Z_s,J_s)=(Z^2,J^2)
\quad&\text{for all $s\geq 1$}.
\endalign
$$
Note that for any $s$ the almost complex structure $J_s$ is compatible
with $\omega$ and we have a smooth family of induced Riemannian metrics $g_s$
on $M$. For a curve $u:\R\to\fl$ we define the energy by
$$
E(u):=\int_\R|u'(s)|^2_{\tilde g_s}ds
=\int_{\R\times S^1}|\ddsu|^2_{g_s}ds dt
=\int_{\R^\times S^1}\omega(\ddsu,J_s\ddsu)dsdt.
$$
$Z$ provides an $s$-dependent vector field on $\fl$, whose flow lines
$u:\R\times S^1\to M$ are the solutions of
$$
D_{Z,J}u:=\ddsu+J_s\ddtu-J_sZ_{s,t}(u(s,t))=0.
\tag{5.1}
$$
As in Proposition~4.1 one shows that for solutions $u:\R\to\fly$ of
\thetag{5.1}, having finite
energy is equivalent to the existence of $x^1\in\pyo$ and $x^2\in\pyt$,
such that $u(s,\cdot)$ converges to $x^1$ \resp $x^2$ when $s\to\infty$
\resp $-\infty$, in the same sense as in Proposition~4.1. For $\tilde
x^1\in\cpyo$ and $\tilde x^2\in\cpyt$ we let $\mot(\tilde x^1,\tilde x^2)$
denote the space of finite energy solutions of \thetag{5.1} connecting $\tilde
x^1$ with $\tilde x^2$. If we want to emphasize the dependence on $(Z,J)$
we will write $\mot(\tilde x^1,\tilde x^2,Z,J)$.
For $u\in\mot(\tilde x^1,\tilde x^2)$ one has an energy estimate
(\cf the proof of Theorem~4.3 in \cite{LO95})
$$
E(u)\leq a^2(\tilde x^2)-a^1(\tilde x^1)+
2\max_{
\Sb
s\in[-1,1],\
t\in S^1
\\x\in M
\endSb
}\Big|\frac{\partial h_{s,t}}
{\partial s}(x)\Big|,
\tag{5.2}
$$
where $i_{Z_{s,t}}\omega=i_{Z^1_t}\omega+dh_{s,t}$ and $h_{s,t}(y(0))=0$.

\definition{Definition 5.1}
Let $(M,\omega)$ be a closed weakly monotone symplectic manifold.
A homotopy $(Z,J)$ between two $\gamma$-regular pairs $(Z^1,J^1)$ and
$(Z^2,J^2)$ is called $\gamma$-regular if the following conditions hold:
\roster
\item
For all $\tilde x^1\in\cpyo$, $\tilde x^2\in\cpyt$ and
$u\in\mot(\tilde x^1,\tilde x^2)$ the linearization of the Cauchy-Riemann
operator $D_{Z,J}$ at $u$ is surjective, \cf \cite{HS95} for definition.
\item
$J:\R\to\Cal J$ is regular in the sense, that
$$
\bigl\{(s,v)\bigm|
\text{$s\in\R$ and $v:S^2\to M$ is $J_s$-holomorphic curve}\bigr\}
$$
is a manifold.
\item
For all $\tilde x^1\in\cpyo$, $\tilde x^2\in\cpyt$ with
$\mu^2(\tilde x^2)-\mu^1(\tilde x^1)\leq 1$, all
$u\in\mot(\tilde x^1,\tilde x^2)$ and all $(s,t)\in\R\times S^1$ on has
$u(s,t)\notin M_0(\infty,J_s)$.
\endroster
\enddefinition

As in Proposition~4.5 one shows

\proclaim{Proposition 5.2}
Let $(M,\omega)$ be a closed weakly monotone symplectic manifold.
\roster
\item
The set of $\gamma$-regular homotopies between two given $\gamma$-regular
pairs is generic in the sense of Baire.
\item
For any $\gamma$-regular homotopy the space
$\mot(\tilde x^1,\tilde x^2)$ is an orientable smooth manifold of dimension
$\mu^2(\tilde x^2)-\mu^1(\tilde x^1)$, for all
$\tilde x^1\in\cpyo$ and $\tilde x^2\in\cpyt$.
\footnote{Notice that this index difference does not depend on the c-structure
above $\gamma$, \cf Remark~3.4.
}
\endroster
\endproclaim

The uniform energy estimate \thetag{5.2} and weak monotonicity gives:

\proclaim{Proposition 5.3}
Let $(M,\omega)$ be a closed weakly monotone symplectic manifold and
$(Z,J)$ a $\gamma$-regular homotopy. Then $\mot(\tilde x^1,\tilde x^2)$ is
a finite set, for all $\tilde x^1\in\cpyo$, $\tilde
x^2\in\cpyt$ with $\mu^2(\tilde x^2)-\mu^1(\tilde x^1)=0$. Moreover
$$
\bigl\{A\in\Gamma\bigm|\phi(A)\leq c, \bar c_1(A)=0,
\mot(\tilde x^1,A\sharp\tilde x^2)\neq\emptyset\bigr\}
$$
is finite for all $c\in\R$ and all $\tilde x^1\in\cpyo$, $\tilde
x^2\in\cpyt$ for which $\mu^2(\tilde x^2)-\mu^1(\tilde x^1)=0$.
\endproclaim

Using coherent orientations one defines
$$
n^{12}(\tilde x^1,\tilde x^2)
:=\sharp\bigl(\mot(\tilde x^1,\tilde x^2)\bigr)\in\Z,
$$
where the points are counted with signs according to their orientation, and
if $\mu^2(\tilde x^2)-\mu^1(\tilde x^1)\neq 0$ we set $n^{12}(\tilde
x^1,\tilde x^2):=0$. Then one has
$$
n^{12}(A\sharp\tilde x^1,A\sharp\tilde x^2)=n^{12}(\tilde x^1,\tilde x^2)
\tag{5.3}
$$
and a gluing argument shows
$$
\sum_{\tilde x\in\cpyo}n^1(\tilde x^1,\tilde x)n^{12}(\tilde x,\tilde x^2)
=
\sum_{\tilde x\in\cpyt}n^{12}(\tilde x^1,\tilde x)n^2(\tilde x,\tilde
x^2),
\tag{5.4}
$$
for all $\tilde x^1\in\cpyo$, $\tilde x^2\in\cpyt$ and all
$A\in\Gamma$. Proposition~5.3 gives
$$
h^{12}(\delta_{\tilde x^1}):=\sum_{\tilde x^2\in\cpyt}n^{12}(\tilde
x^1,\tilde x^2)\delta_{\tilde x^2}\in C_F^*(Z^2,J^2).
$$
From \thetag{5.3} we see that $h^{12}(\delta_A*\delta_{\tilde
x^1})=\delta_A*h^{12}(\delta_{\tilde x^1})$, and so $h^{12}$
extends uniquely to a $\Lambda$-linear map
$$
h^{12}:C_F^*(Z^1,J^1)\to C_F^*(Z^2,J^2)
$$
of degree 0. Moreover \thetag{5.4} immediately gives
$\partial^2\circ h^{12}=h^{12}\circ\partial^1$, so $h^{12}$ is a chain map
and induces
$$
h^{12}:H_F^*(Z^1,J^1)\to H_F^*(Z^2,J^2).
$$

Suppose one has two $\gamma$-regular homotopies $(Z_0,J_0)$ and
$(Z_1,J_1)$ connecting the same $\gamma$-regular pairs $(Z^1,J^1)$ and
$(Z^2,J^2)$. Then a homotopy of homotopies between these two
$\gamma$-regular homotopies is a pair $(Z,J)$, where
$Z\in C^\infty([0,1]\times\R,\vf^\beta)$ and
$J\in C^\infty([0,1]\times\R,\Cal J)$, such that
$$
\align
(Z_{i,s},J_{i,s})&=(Z_i,J_i)\quad
\text{for all $s\in\R$ and $i=0,1$,}
\\
(Z_{\lambda,s},J_{\lambda,s})&=(Z^1,J^1)\quad
\text{for all $\lambda\in[0,1]$, $s\leq -1$ and}
\\
(Z_{\lambda,s},J_{\lambda,s})&=(Z^2,J^2)\quad
\text{for all $\lambda\in[0,1]$, $s\geq 1$.}
\endalign
$$
For $\tilde x^1\in\cpyo$ and $\tilde x^2\in\cpyt$ define
$$
{\Cal H}(\tilde x^1,\tilde x^2):=\bigl\{
(\lambda,u)\bigm|\lambda\in[0,1], u\in\mot(\tilde x^1,\tilde
x^2,Z_\lambda,J_\lambda)\bigr\}.
$$

\definition{Definition~5.4}
A homotopy $(Z,J)$ between two $\gamma$-regular homotopies connecting the same
two $\gamma$-regular pairs as above is called $\gamma$-regular if the following
holds:
\roster
\item
For all $\tilde x^1\in\cpyo$ and $\tilde x^2\in\cpyt$ the space
$\Cal H(\tilde x^1,\tilde x^2)$ is a manifold of dimension
$\mu^2(\tilde x^2)-\mu^1(\tilde x^1)+1$.
\item
$J:[0,1]\times\R\to\Cal J$ is regular in the sense that
$$
\bigl\{(\lambda,s,v)\bigm|
\text{$\lambda\in [0,1]$, $s\in\R$ and $v:S^2\to M$ is a
$J_{\lambda,s}$-holomorphic curve}\bigr\}
$$
is a manifold.
\item
For all $\tilde x^1\in\cpyo$, $\tilde x^2\in\cpyt$ with
$\mu^2(\tilde x^2)-\mu^1(\tilde x^1)+1\leq 1$, all
$(\lambda,u)\in\Cal H(\tilde x^1,\tilde x^2)$ and all $(s,t)\in\R\times
S^1$ one has $u(s,t)\notin M_0(\infty,J_{\lambda,s,t})$.
\endroster
\enddefinition

\proclaim{Proposition 5.5}
The set of $\gamma$-regular homotopies between two $\gamma$-regular homotopies
connecting the same $\gamma$-regular pairs is generic in the sense of Baire.
\endproclaim

\proclaim{Proposition 5.6}
For a $\gamma$-regular homotopy of homotopies $(Z,J)$ the manifold
${\Cal H}(\tilde x^1,\tilde x^2)$ is a finite set, for all
$\tilde x^1\in\cpyo$ and $\tilde x^2\in\cpyt$ with
$\mu^2(\tilde x^2)-\mu^1(\tilde x^1)+1=0$. Moreover
$$
\bigl\{A\in\Gamma\bigm|\phi(A)\leq c, \bar c_1(A)=0,{\Cal H}(\tilde
x^1,A\sharp \tilde x^2)\neq\emptyset\bigr\}
$$
is finite for all $c\in\R$ and all $\tilde x^1\in\cpyo$, $\tilde
x^2\in\cpyt$ with $\mu^2(\tilde x^2)-\mu^1(\tilde x^1)+1=0$.
\endproclaim

Using coherent orientations one defines
$$
m(\tilde x^1,\tilde x^2)
:=\sharp\bigl({\Cal H}(\tilde x^1,\tilde x^2)\bigr)\in\Z
$$
for $\mu^2(\tilde x^2)-\mu^1(\tilde x^1)+1=0$
and sets $m(\tilde x^1,\tilde x^2)=0$ otherwise. Then one has
$$
m(A\sharp\tilde x^1,A\sharp\tilde x^2)=m(\tilde x^1,\tilde x^2)
\tag{5.5}
$$
and from a gluing argument one gets
$$
n^{12}_1(\tilde x^1,\tilde x^2)-n^{12}_0(\tilde x^1,\tilde x^2)=
\sum_{\tilde x\in\cpyt}\Big(m(\tilde x^1,\tilde x)n^2(\tilde x,\tilde
x^2)
+n^1(\tilde x^1,\tilde x)m(\tilde x,\tilde x^2)\Big).
\tag{5.6}
$$
Because of Proposition~5.6 and \thetag{5.5}
$$
H(\delta_{\tilde x^1}):=\sum_{\tilde x^2\in\cpyt}m(\tilde x^1,\tilde
x^2)\delta_{\tilde x^2}
$$
extends uniquely to a $\Lambda$-linear map $H:C^*_F(Z^1,J^1)\to
C_F^{*-1}(Z^2,J^2)$ of degree $-1$. Equation \thetag{5.6} shows that $H$ is a
chain homotopy, \ie $h^{12}_1-h^{12}_0=\partial^2 H+H\partial^1$. So
$$
h^{12}_1=h^{12}_0:H_F^*(Z^1,J^1)\to H_F^*(Z^2,J^2).
\tag{5.7}
$$

Using the constant homotopy, which is regular in the sense
of Definition~5.1 one immediately gets $h^{11}=\id$, even on chain level for
this special homotopy.
Concatenating two homotopies and using a gluing argument
one shows $h^{23}\circ h^{12}=h^{13}$, again on chain level if the
ambiguity in the concatenation is chosen large enough, \cf Lemma~6.4 in
\cite{SZ92}.
It follows immediately that \thetag{5.7} are canonic isomorphisms and one
gets, \cf \cite{F89}, \cite{SZ92}, \cite{HS95} and \cite{LO95}:

\proclaim{Theorem 5.7}
Let $(M,\omega)$ be a closed connected weakly monotone symplectic manifold,
$\beta\in H^1(M;\R)$, $[\tilde y]$ a c-structure above $\gamma$, and $R$ a
commutative ring
with unit. Any $\gamma$-regular pair $(Z,J)\in\vfbyreg\times\Jreg$ and
coherent orientations $\Cal O$ provide a $\Z_{2N}$-graded cochain complex
$C^*_F$ of free $\Lambda$-modules. Any $\gamma$-regular homotopy between such
pairs induces a morphism of cochain complexes and an isomorphism in cohomology.
Different homotopies induce chain homotopic morphisms and therefore
the cohomologies of all these cochain complexes are canonically isomorphic
and will be denoted by
$$
H_F^*(M,\omega,\beta,[\tilde y],R).
$$
Here $\Lambda$ is the Novikov ring associated to the group
$\Gamma$ and the weighting homomorphism
$\phi:\Gamma\to\R$, and $N$ is the minimal Chern number defined by $\bar
c_1\big(H_1(\fly;\Z)\big)=N\Z$. The associated $\Z_2$-graded module
depends only on $(M,\omega,\beta,\gamma,R)$.
\endproclaim


Let $g$ be any Riemannian metric on the compact manifold $M$. For 
$a\in\Omega^*(M)$ define the supremum norm 
$||a||_0:=\sup_{x\in M}|a_x|$ and consider the norm $|\cdot|$ on
$H^*(M,\R)$ defined by $|\alpha|:=||a||_0$, where $a$ is the harmonic
representative of $\alpha\in H^*(M;\R)$. Let $r(g)>0$
denote the injectivity radius of $(M,g)$. If $(M,\omega)$ is symplectic we
set $\varepsilon(g,\omega):=r(g)/||\sharp_\omega||_0>0$, where
$\sharp_\omega:T^*M\to TM$ denotes the vector bundle isomorphism induced by 
$\omega$. Suppose $\beta\in H^1(M;\R)$ with $|\beta|\leq\varepsilon(g,\omega)$.
Consider the time independent symplectic vector field $Z:=\sharp_\omega b$, 
where $b\in\Omega^1(M)$ is the harmonic representative of $\beta$. Then 
$||Z||_0\leq r(g)$, hence any $1$-periodic trajectory of $Z$ must be 
contractible.

\proclaim{Theorem 5.8}
Let $(M,\omega)$ be a closed connected weakly monotone symplectic
manifold, $\gamma$ be a homotopy class of closed curves in $M$ and suppose
that at least one of the following conditions is satisfied:
\roster
\item
$\beta(\gamma)\neq 0$ or
\item
$\gamma$ non-trivial and $|\beta|\leq\varepsilon(\omega,g)$ for any Riemannian metric $g$ on $M$.
\endroster
Then
$$
H_F^*(M,\omega,\beta,[\tilde y],R)=0,
$$
for all c-structures $[\tilde y]$ above $\gamma$ and all commutative rings 
with unit $R$.
\endproclaim

%

\demo{Proof}
Any of the conditions immediately implies, that there exists a time
independent symplectic vector field $Z$, $[i_Z\omega]=\beta$, which can not 
have $1$-periodic trajectories in the homotopy class of
$\gamma$. Thus, for any $J\in\Jreg$ the pair $(Z,J)$ is
$\gamma$-regular in the sense of Definition~4.4, and the Floer complex
$C^*_F(M,\omega,[\tilde y],Z,J,\Cal O,R)$ vanishes. In view of Theorem~5.7 
we get $H_F^*(M,\omega,\beta,[\tilde y],R)=0$.
\qed
\enddemo

\remark{Remark}
It is possible to construct $\beta$ and $\gamma$ both non-trivial but
$\beta(\gamma)=0$, so that the Floer cohomology is non-trivial.
\endremark

\proclaim{Corollary 5.9}
Let $(M,\omega)$ be a closed connected weakly monotone symplectic  manifold,
and suppose $\gamma$ is a homotopy class of closed curves in
$M$. Moreover let $Z$ be a time dependent
symplectic vector field, such that $[i_{Z_t}\omega]$ is a fixed class in
$\beta\in H^1(M;\R)$ and such that all 
$1$-periodic trajectories in $\py$ are non-degenerate. If 
$H^*_F(M,\omega,\beta,\gamma,\R)=0$ then
$$
|\pyeven|=|\pyodd|,
$$
where $|\pyeven|$ \resp $|\pyodd|$ denotes the number of elements
of even \resp odd Conley-Zehnder index in $\py$.
\endproclaim

\demo{Proof}
By Proposition~4.5\therosteritem2 we may assume that there exists
$J\in\Cal J$, such that $(Z,J)$ is a $\gamma$-regular pair.
Choosing coherent orientations $\Cal O$ we get a well defined Floer 
complex $C^*_F(M,\omega,\beta,\gamma,Z,J,\Cal O,\R)$, whose 
cohomology vanishes by Theorem~5.7. Since
$C^*_F(M,\omega,\beta,\gamma,Z,J,\Cal O,\R)$ is a free
acyclic complex its Euler characteristic must vanish, \ie $|\pyeven|=|\pyodd|$.
\qed
\enddemo

From Corollary~5.9 we immediately get

\proclaim{Corollary 5.10}
Let $(M,\omega)$ be closed connected weakly monotone, and let $Z$ be any
$1$-periodic time dependent symplectic  vector field, such that
$[i_{Z_t}\omega]$ is a fixed class $\beta\in H^1(M;\R)$. If $Z$ has a
non-degenerate $1$-periodic trajectory in the homotopy class $\gamma$ and 
$H^*_F(M,\omega,\beta,\gamma,\R)=0$,
then it must have another $1$-periodic solution in the 
same homotopy class, which might of course be degenerate.
\endproclaim

The last two corollaries prove Theorem~1.1 stated in the introduction.

\remark{Remark}
Both corollaries remain true for non-trivial $\gamma$ without the 
hypothesis $H^*_F(M,\omega,\beta,\gamma,\R)=0$. They were also known to
K.~Ono.
\endremark


\remark{Remark~5.11}
Theorem~5.7 and 5.8 remain true for any closed symplectic manifold,
and any
pair $(Z,J)\in\vfbyreg\times\Cal J$
and therefore Corollary~5.9 and 5.10 remain true for any closed symplectic manifold.
\endremark

\head
6. Symplectic torsion
\endhead

In this section we introduce a new invariant, the symplectic torsion. It is
associated with a closed symplectic manifold $(M,\omega)$, a homotopy class
of closed curves $\gamma$ and a pair $(Z,J)\in\vfbyreg\times\Cal J$. The
invariant takes value in the Abelian group $\Wh(\Lambda)$ defined below.

Since the theory of Floer complex was considered only for weakly monotone manifolds
and for $\gamma$-regular pairs, these assumptions will be understood here too, but
they are not necessary.

We recall some basic properties and definitions of the Milnor torsion. Let $A$
be (for simplicity) a commutative ring with unit. We consider the Abelian groups
$$
K_1(A):=\frac{\GL_\infty(A)}{[\GL_\infty(A),\GL_\infty(A)]}
\quad\text{and}\quad
\bar K_1(A):=K_1(A)/(\pm 1),
$$
where $(\pm 1)\in K_1(A)$ are the elements represented by
$(\pm 1)\in\GL_1(A)$. Recall that one has a
homomorphism $\det:\bar K_1(A)\to U(A)/\{\pm 1\}$,
where $U(A)$ denotes the units of $A$, and if $A$ has no zero divisors
this is an isomorphism. We will write $\bar K_1(A)$ additively.

Let $M$ be a finitely generated free $A$-module
and suppose $b=\{b_1,\dotsc,b_r\}$
and $c=\{c_1,\dotsc,c_r\}$ are two bases of $M$. Then one has
$b_i=\sum_{j=1}^ka_{ij}c_j$ for some $a_{ij}\in A$ and we will write
$[b/c]\in\bar K_1(A)$ for the element represented by the $r\times r$-matrix
$(a_{ij})$. Note that this does not depend on the ordering of the bases, and
$[b/d]=[b/c]+[c/d]$ for bases $b$, $c$ and $d$.

Now let $M^*$ be a free finitely generated $\Z_2$-graded $A$-module. A graded
base for $M^*$ is $b^*=(b^0,b^1)$, where $b^i$ is a base of $M^i$, $i\in\Z_2$.
For two graded bases $b^*$ and $c^*$ we define
$$
[b^*/c^*]:=[b^0/c^0]-[b^1/c^1]\in\bar K_1(A).
$$
We call two graded bases $b^*$ and $c^*$ equivalent if $[b^*/c^*]=0$
and will write $\bold b^*$ for the equivalence class. Note that
$[\bbs/\bold c^*]\in\bar K_1(A)$ is well defined, and one has
$$
[\bbs/\bds]=[\bbs/\bcs]+[\bcs/\bds]
\tag{6.1}
$$
for three equivalence classes of graded bases $\bbs$, $\bcs$ and $\bds$ on
$M^*$. Suppose one has a short exact sequence of free finitely generated
$\Z_2$-graded $A$-modules
$$
0\to M^*_1\to M^*_2\to M^*_3\to 0
$$
and suppose $\bbs_1$ \resp $\bbs_3$ are equivalence classes of graded
bases on $M^*_1$ \resp $M^*_3$. Choosing representing graded bases
$b^*_1$, $b^*_3$ and lifting $b^*_3$ to $M^*_2$ one obtains a graded
base $b^*_1b^*_3$ of $M^*_2$ and a well defined equivalence class of
graded bases $\bbs_1\bbs_3$ on $M^*_2$. If $\tilde\bbs_1$ and $\tilde\bbs_3$
are other equivalence classes of graded bases on $M^*_1$ and
$M^*_3$, one has
$$
[\tilde\bbs_1\tilde\bbs_3/\bbs_1\bbs_3]
=[\tilde\bbs_1/\bbs_1]+[\tilde\bbs_3/\bbs_3].
\tag{6.2}
$$

Next consider a free finitely generated $\Z_2$-graded chain complex of
$A$-modules $(C^*,\partial^*)$. It gives rise to three $\Z_2$-graded
$A$-modules, namely the cycles $Z^*:=\ker(\partial^*)$, the boundaries
$B^*:=\img(\partial^{*-1})$ and the homology $H^*=Z^*/B^*$. Suppose
$H^*$ and $B^*$ are also free, and $C^*$ and $H^*$ are equipped with
equivalence classes of graded bases $\bcs$ and $\bhs$, respectively. Choose
an equivalence class of graded bases $\bbs$ of $B^*$. The short exact sequences
of $\Z_2$-graded $A$-modules
$$
0\to B^*\to Z^*\to H^*\to 0
\quad\text{and}\quad
0\to Z^*\to C^*@>{\partial^*}>>B^{*+1}\to 0
$$
give rise to an equivalence class of graded bases $\bbs\bhs$ on $Z^*$ and
$(\bbs\bhs)\bbss$ on $C^*$.
\footnote{
If $B^*$ is a $\Z_2$-graded module $B^{*+1}$ denotes the same module with
shifted grading. So $\partial^*:C^*\to B^{*+1}$ is of degree $0$. If $b^*$
is a graded base for $B^*$, then $b^{*+1}$ will denote the corresponding
graded base of $B^{*+1}$, and similarly for equivalence classes of
graded bases.}
The Milnor torsion is now defined by, \cf~\cite{Mi66},
$$
\tau(C^*,\partial^*,\bcs,\bhs):=[(\bbs\bhs)\bbss/\bcs]\in\bar K_1(A).
$$
A straight forward calculation using \thetag{6.1}, \thetag{6.2} and
$[\bbs/\tilde\bbs]=-[\bbss/\tilde\bbss]$ shows that it does not depend on the
choice of $\bbs$, and if $\tilde\bcs$ and $\tilde\bhs$ are other
equivalence classes of graded bases one gets
$$
\tau(C^*,\partial^*,\tilde\bcs,\tilde\bhs)-\tau(C^*,\partial^*,\bcs,\bhs)
=[\tilde\bhs/\bhs]-[\tilde\bcs/\bcs].
\tag{6.3}
$$
Suppose one has a short exact sequence of free $\Z_2$-graded acyclic chain
complexes $0\to C_1^*\to C_2^*\to C_3^*\to 0$, such that $B_1^*$ and
$B^*_3$ are free with equivalence classes of graded bases $\bcs_1$ and
$\bcs_3$. Then $\bcs_1\bcs_3$ is an equivalence class of graded bases of
$C^*_2$ and
$$
\tau(C^*_2,\partial^*_2,\bcs_1\bcs_3,\emptyset)
=\tau(C^*_1,\partial^*_1,\bcs_1,\emptyset)
+\tau(C^*_3,\partial^*_3,\bcs_3,\emptyset).
\tag{6.4}
$$
To see \thetag{6.4} recall that one has a short exact sequence
$0\to B_1^*\to B_2^*\to B_3^*$, choose equivalence classes of graded bases
$\bbs_1$ and $\bbs_3$ of $B_1^*$ and $B_3^*$ and equip $B_2^*$ with
$\bbs_1\bbs_3$. By definition we have
$\tau(C^*_2,\partial^*_2,\bcs_1\bcs_3,\emptyset)=
[(\bbs_1\bbs_3)(\bbss_1\bbss_3)/\bcs_1\bcs_3]$. It is easy to see, that
$[(\bbs_1\bbs_3)(\bbss_1\bbss_3)/(\bbs_1\bbss_1)(\bbs_3\bbss_3)]=0$.
Together with \thetag{6.1} and \thetag{6.2} this immediately gives
\thetag{6.4}.

Given a chain mapping $f^*:(C_1^*,\partial^*_1)\to(C_2^*,\partial^*_2)$
between two $\Z_2$-graded complexes, we consider the `mapping cone'
$(C^*_f,\partial_f^*)$, a $\Z_2$-graded chain complex defined by
$$
C^*_f:=C^*_2\oplus C^{*+1}_1,
\quad
\partial_f^*:=
\left(
  \smallmatrix
    \partial^*_2 & (-1)^{*+1}f^{*+1}\\
    0&\partial_1^{*+1}
  \endsmallmatrix
\right).
$$
Note, that one has a short exact sequence
$0\to C^*_2\to C^*_f\to C^{*+1}_1\to 0$ of $\Z_2$-graded complexes.
We assume that $C^*_1$, $C^*_2$, $B^*_1$ and $Z^*_2$ are free and that $f^*$
induces an isomorphism in homology. Then $C^*_f$ is free, acyclic and
$Z_f^*=B^*_f$ is free, too. Given equivalence classes of graded bases $\bcs_1$
and $\bcs_2$ of $C^*_1$ and $C^*_2$ one gets an equivalence class of graded
bases $\bcs_2\bcss_1$ on $C^*_f$ and defines the relative torsion by
$$
\tau(f^*,\bcs_1,\bcs_2)
:=\tau(C^*_f,\partial^*_f,\bcs_2\bcss_1,\emptyset)\in\bar K_1(A).
$$
From \thetag{6.3} one immediately gets
$$
\tau(f^*,\tilde\bcs_1,\tilde\bcs_2)-\tau(f^*,\bcs_1,\bcs_2)
=[\tilde\bcs_1/\bcs_1]-[\tilde\bcs_2/\bcs_2].
\tag{6.5}
$$
If $f$ and $g$ are chain homotopic, \ie there exists $H:C^*_1\to
C^{*-1}_2$ satisfying $f^*-g^*=\partial_2^{*-1}H^*+H^{*+1}\partial^*_1$,
then one obtains an isomorphism of chain complexes
$$
\psi_H^*:C^*_f=C^*_2\oplus C^{*+1}_1\to C^*_2\oplus C^{*+1}_1=C^*_g,
\quad
\psi^*_H=
\left(\smallmatrix
\id_{C^*_2} & (-1)^{*+1}H^{*+1} \\
0 & \id_{C_1^{*+1}}
\endsmallmatrix\right).
$$
Obviously $\psi^*_H(\bcs_2\bcss_1)=\bcs_2\bcss_1$ and hence
$$
\tau(f^*,\bcs_1,\bcs_2)=\tau(g^*,\bcs_1,\bcs_2)
\tag{6.6}
$$
for chain homotopic maps $f$ and $g$.

Suppose we have $f^*:(C_1^*,\partial^*_1)\to(C_2^*,\partial^*_2)$ and
$g^*:(C_2^*,\partial^*_2)\to(C_3^*,\partial^*_3)$ as above and
equivalence classes of graded bases $\bcs_1$, $\bcs_2$ and $\bcs_3$.
Consider the chain mappings
$$
\psi^*:=\left(\smallmatrix g & 0 \\ 0 & \id_{C^*_1}\endsmallmatrix\right):
C^*_f\to C^*_{gf}
\quad\text{and}\quad
\phi^*:=\left(\smallmatrix gf & 0 \\ 0 & f\endsmallmatrix\right):
C^*_{\id_{C_1^*}}\to C^*_g
$$
These give rise to short exact sequences
$0\to C^*_{gf}\to C^*_\psi\to C^{*+1}_f\to 0$ and
$0\to C^*_g\to C^*_\phi\to C^{*+1}_{\id_{C_1^*}}\to 0$.
Interchanging the second and third factor defines an isomorphism from
$C^*_\psi$ to $C^*_\phi$ which sends $\bcs_3\bcss_1\bcss_2\bcsss_1$ to
$\bcs_3\bcss_2\bcss_1\bcsss_1$. So
$\tau(\psi^*,\bcs_2\bcss_1,\bcs_3\bcss_1)
=\tau(\phi^*,\bcs_1\bcss_1,\bcss_3\bcs_2)$. Together with \thetag{6.4}
this yields
$$
\tau(g^*\circ f^*,\bcs_1,\bcs_3)
=\tau(f^*,\bcs_1,\bcs_2)+\tau(g^*,\bcs_2,\bcs_3),
\tag{6.7}
$$
since obviously $\tau(\id_{C^*_1},\bcs_1,\bcs_1)=0$.

If $f^*:(C^*_1,\partial^*_1)\to(C^*_2,\partial^*_2)$ is a chain mapping
between acyclic $\Z_2$-graded complexes, such that $B^*_1$, $B^*_2$ are free,
then the short exact sequence $0\to C^*_2\to C^*_f\to C^{*+1}_1\to 0$ and
\thetag{6.4} yield
$$
\tau(f^*,\bcs_1,\bcs_2)=
\tau(C^*_2,\partial^*_2,\bcs_2,\emptyset)
-\tau(C^*_1,\partial^*_1,\bcs_1,\emptyset).
\tag{6.8}
$$

We will apply the previous definitions to the Floer complex defined in
section~4. Recall that the Floer complex was a free $\Z_2$-graded cochain
complex over the Novikov ring $\Lambda$ associated to
$(M,\omega,\gamma,Z,J,\Cal O,R)$.
We will assume that $R$ is a principal ideal domain. Then $\Lambda$ has
no zero divisors and the subring $\Lambda_0$ is a principal ideal domain
too, see section~2. Note that via $\Gamma\subseteq
U(\Lambda)=\GL_1(\Lambda)$ \resp $\Gamma_0\subseteq
U(\Lambda_0)=\GL_1(\Lambda_0)$, $\Gamma$ \resp $\Gamma_0$ becomes a
subgroup of $K_1(\Lambda)$ \resp $K_1(\Lambda_0)$.
We define the Whitehead groups by
$$
\Wh(\Lambda):=\frac{K_1(\Lambda)}{\{\pm\Gamma\}}
=\frac{\bar K_1(\Lambda)}{\Gamma}
\quad\text{and}\quad
\Wh(\Lambda_0):=\frac{K_1(\Lambda_0)}{\{\pm\Gamma_0\}}
=\frac{\bar K_1(\Lambda_0)}{\Gamma_0}.
$$
Note that since $\Lambda$ and $\Lambda_0$ have no zero divisors we have
isomorphisms
$$
\det:\Wh(\Lambda)\cong U(\Lambda)/\{\pm\Gamma\}
\quad\text{and}\quad
\det:\Wh(\Lambda_0)\cong U(\Lambda_0)/\{\pm\Gamma_0\}.
$$
So the inclusion $i:\Lambda_0\to\Lambda$ induces an injective homomorphism
$\Wh(i):\Wh(\Lambda_0)\to\Wh(\Lambda)$.

\subhead
The relative torsion $\tau\bigl(\gamma,(Z_1,J_1),(Z_2, J_2)\bigr)$
\endsubhead
To the data $(M,\omega,\gamma,Z,J,\Cal O,R)$, where $(Z,J)$ is
$\gamma$-regular, we have in (section~4) associated the Floer complex
$C^*_F$, a free $\Z_2$-graded cochain complex of $\Lambda$-modules.
Choosing a lift $\tilde x\in\cpy$ for
every $x\in\py$ defines a base $\{\delta_{\tilde x}\mid x\in\py\}$ of
$C^*_F$. One can
choose them, such that the matrix expression of $\partial^0:C^0_F\to C^1_F$
has entries in $\Lambda_0$. So there exist free $\Lambda_0$-modules
$\bar C^0_F$, $C^1_F$ and
$\bar\partial^0:\bar C^0_F\to\bar C^1_F$, such that
$$
C^0_F=\Lambda\otimes_{\Lambda_0}\bar C^0_F,
\quad
C^0_F=\Lambda\otimes_{\Lambda_0}\bar C^1_F
\quad\text{and}\quad
\partial^0=\Lambda\otimes_{\Lambda_0}\bar\partial^0.
$$
Since $\Lambda_0$ is a
principal ideal domain we see that
$Z^0_F=\Lambda\otimes_{\Lambda_0}\bar Z^0_F$ and
$B^1_F=\Lambda\otimes_{\Lambda_0}\bar B^1_F$ are free $\Lambda$-modules.
Similarly, but one has to choose the lifts differently, one shows that
$Z^1_F$ and $B^0_F$ are free $\Lambda$-modules.

For any two $\gamma$-regular pairs $(Z_1,J_1)$ and $(Z_2,J_2)$ with
$\beta_1=\beta_2\in H^1(M;\R)$, we have (in section~5) constructed
chain mappings $h^{12}:C^*_{1,F}\to C^*_{2,F}$ and all of them were chain
homotopic, \cf Theorem~5.7. For every $x^1\in\pyo$ and $x^2\in\pyt$ choose
lifts $\tilde x^1\in\cpyo$ and $\tilde x^2\in\cpyt$. So we get bases
$\bcs_1:=\{\delta_{\tilde x^1}\mid x^1\in\pyo\}$ of $C_F^*(Z_1,J_1)$ and
$\bcs_2:=\{\delta_{\tilde x^2}\mid x^2\in\pyt\}$ of $C_F^*(Z_2,J_2)$.
For two $\gamma$-regular pairs $(Z_1,J_1)$ and $(Z_2,J_2)$ with
$\beta_1=\beta_2\in H^1(M;\R)$ we define the relative torsion by
$$
\tau\bigl(\gamma,(Z_1,J_1),(Z_2,J_2)\bigr):=\tau(h^{12},\bcs_1,\bcs_2)
\in\Wh(\Lambda).
$$
By \thetag{6.6} this does not depend on $h^{12}$ and by \thetag{6.5} it does
not depend on the lifts $\tilde x^1$ and $\tilde x^2$, for different choices
give rise to different bases $\tilde\bcs_1$, $\tilde\bcs_2$ which obviously
satisfy $[\tilde\bcs_1,\bcs_1]\in\Gamma$ and $[\tilde\bcs_2,\bcs_2]\in\Gamma$.
From \thetag{6.7} and the remarks right before Theorem~5.7 we get
$$
\tau\bigl(\gamma,(Z_1,J_1),(Z_3,J_3)\bigr)
=\tau\bigl(\gamma,(Z_1,J_1),(Z_2,J_2)\bigr)
+\tau\bigl(\gamma,(Z_2,J_2),(Z_3,J_3)\bigr),
$$
for every $\gamma$-regular pairs $(Z_1,J_1)$, $(Z_2,J_2)$ and $(Z_3,J_3)$,
for which $\beta_1=\beta_2=\beta_3\in H^1(M;\R)$.

\remark{Remark 6.1}
Since $\phi:\Gamma_0\to\R$ is injective every
$\lambda\in\Lambda_0$ has a well defined `leading term'
$\lambda(g)\in R$, where $g\in\Gamma_0$ is the unique element in $\Gamma_0$,
such that $\lambda(g)\neq 0$ and $\phi(g)$ is minimal with that property.
If $\lambda=0$ we define it to be $0\in R$. In \cite{HS95} it is shown
that $\lambda\in\Lambda_0$ is invertible if and only if this
leading term is invertible in $R$. So one obtains a homomorphism
$$
\lt:U(\Lambda_0)/\{\pm\Gamma_0\}\to U(R)/\{\pm 1\}.
$$
\endremark

A little inspection of the Floer complex and Remark~6.1 above permits to show

\proclaim{Proposition 6.2}
\roster
\item
The torsion $\tau\bigl(\gamma,(Z_1,J_1),(Z_2,J_2)\bigr)$ lies in the image
of the injective homomorphism $\Wh(\Lambda_0)\to\Wh(\Lambda)$.
\item
The homomorphism
$\Wh(\Lambda_0)@>\det>>U(\Lambda_0)/\{\pm\Gamma_0\}@>\lt>>U(R)/\{\pm 1\}$
sends the torsion $\tau\bigl(\gamma,(Z_1,J_1),(Z_2,J_2)\bigr)$ to $1$.
\endroster
\endproclaim

\subhead
The symplectic torsion for $H_F^*(M,\omega,\beta,\gamma,R)=0$
\endsubhead
Since the Floer complex is acyclic we can consider the
Milnor torsion associated to the Floer complex generated by the
$\gamma$-regular pair $(Z,J)$ and coherent orientations $\Cal O$.
Together with the choice of the lifts $\tilde x$ these data will lead to an
element in $\Wh(\Lambda)$ denoted by
$$
\tau\bigl(\gamma,(Z,J)\bigr)\in\Wh(\Lambda),
$$
independent of orientations, lifts and the c-structure above $\gamma$. The
relation to the relative torsion is given by, see \thetag{6.8},
$$
\tau\bigl(\gamma,(Z_1,J_1),(Z_2,J_2)\bigr)
=\tau\bigl(\gamma,(Z_2,J_2)\bigr)-\tau\bigl(\gamma,(Z_1,J_1)\bigr),
$$
for all $\gamma$-regular pairs with $\beta_1=\beta_2\in H^1(M;\R)$.
If $(Z',J')$ is a $\gamma$-regular pair whose vector field $Z'$ has no
$1$-periodic trajectories in the class $\gamma$, then
$\tau\bigl(\gamma,(Z',J')\bigr)=0$ and
$$
\tau\bigl(\gamma,(Z,J)\bigr)=\tau(\gamma,(Z',J'),(Z,J)),
$$
for every $\gamma$-regular pair $(Z,J)$, with $\beta=\beta'\in H^1(M;\R)$.

In section~8 below we will show that this symplectic torsion is in general
non-trivial. The non-triviality of this torsion assures the existence of
closed trajectories of $Z$ in the class $\gamma$, and the non-existence of
such trajectories implies the vanishing of this torsion.

\subhead
The symplectic torsion for $\gamma$ trivial
\endsubhead
Denote by $\gamma_0$ the trivial class, \ie the class of contractible loops.

\definition{Definition 6.3}
Suppose $(Z,J)\in\Cal Y^\beta_{\text{$\gamma_0$-reg}}\times\Jreg$, set
$\alpha:=i_Z\omega$ which is a closed one form and let $g$ be the
Riemannian metric induced by $J$ and $\omega$, \cf section~4.
The pair $(Z,J)$ is called special if the following condition holds:
\roster
\item
The pair $(\alpha,g)$ is Morse-Smale
\footnote{The pair $(\alpha,g)$ consisting of a Riemannian metric $g$ and a
closed one form is Morse-Smale if all critical points of $\alpha$ are
non-degenerate and for any two such critical points $x$ and $y$ the
stable manifold of $x$ and the unstable manifold of $y$ intersect
transversally, \cf \cite{BH01}.}
and the unstable sets of the critical points with respect to the vector
field $\grad_g\alpha$ provide a cell structure of $M$.
\endroster
For the reader unfamiliar with the Morse-Smale pairs \therosteritem1
can be replaced by
\roster
\item"(1')"
There exists a smooth triangulation of $M$ so that the unstable sets
associated with the critical points of $\alpha$ with respect to the vector
field $\grad_g\alpha$ identify to the open simplexes of the triangulation.
\endroster
\enddefinition

\proclaim{Lemma 6.4}
For all $\beta\in H^1(M;\R)$ there exist special pairs.
\endproclaim

As in \cite{FO99} one can define a Floer complex every pair
$(Z,J)\in\vfboreg\times\Jreg$ and Theorem~5.7 remains true, \cf
Remark~5.11. So the relative torsion is also well defined for two such
pairs and Proposition~6.2 still holds. In section~7 we will prove

\proclaim{Proposition 6.5}
Given two special pairs $(Z_1,J_1)$ and $(Z_2,J_2)$ the
relative torsion $\tau\bigl(\gamma_0,(Z_1,J_1),(Z_2,J_2)\bigr)$ vanishes.
\endproclaim

Now given an arbitrary $\gamma_0$-regular pair $(Z,J)$, we define
$$
\tau\bigl(\gamma_0,(Z,J)\bigr):=\tau\bigl(\gamma,(Z',J'),(Z,J)\bigr),
$$
where $(Z',J')$ is any special $\gamma_0$-regular pair. Proposition~6.5
shows, that this definition is independent of the special pair $(Z',J')$.
Proposition~6.2 shows, that the torsion $\tau\bigl(\gamma_0,(Z,J)\bigr)$ is
actually an element in $\Wh(\Lambda_0)$ with vanishing leading term.

As indicated in introduction, if $Z$ is time independent and $\gamma$ is trivial,
then
this torsion can be calculated. It also has  interesting geometric interpretations.
This will be presented in a future paper.

\head
7. Proof of Proposition~6.5
\endhead

Let us first recall a few elements of Novikov theory, \cf~\cite{BH01}. Let
$M$ be closed oriented manifold and let $\alpha$ be a Morse form, \ie $\alpha$
is closed and all critical points (zeros) are non-degenerated. Let
$\beta\in H^1(M;\R)$ denote the cohomology class represented by $\alpha$
and suppose we have a connected principal Abelian covering $\pi:\tilde M\to M$,
such that $\pi^*\beta=0\in H^1(\tilde M;\R)$. Let $\Delta$ denote the
structure group and note, that $\beta$ induces a homomorphism
$\beta:\Delta\to\R$. Choose a commutative ring with unit $R$ and let
$\Cal N:=\Cal N(\Delta,\beta,R)$ denote the corresponding Novikov ring. The
underlying $\Cal N$-module of the Morse-Novikov complex $C^*_N$ is the set
of functions from $\tilde C:=\pi^{-1}(C)$ to $R$ which satisfy a Novikov
condition as in section~3. Here $C$ denotes the critical points of
$\alpha$. $C^*_N$ is a free $\Cal N$-module of rank equal to the
cardinality of $C$ and it is $\Z$-graded by the Morse index.

Now choose a Riemannian metric $g$ on $M$ and assume that $\grad_g\alpha$
is Morse-Smale. Finally choose orientations $\Cal O$
for the unstable manifolds. Counting the trajectories
of $\pi^*(\grad_g\alpha)$ connecting points of index difference $1$ in
$\tilde C$ with appropriate sings defines an $\Cal N$-linear differential
$\partial$ of degree $1$ on $C^*_N$. This complex is called
Morse-Novikov complex and we will denote it for simplicity by
$C^*_N(\alpha,g)$, although it also depends on the orientations $\Cal O$,
the covering $\tilde M$ and the ring $R$.

Given two Morse-Smale pairs $(\alpha_1,g_1)$ and $(\alpha_2,g_2)$ a
homotopy between them is a smooth family of $1$-forms $\alpha_t$
which all represent the same cohomology class $\beta$ and
a smooth family of Riemannian metrics $g_t$, such that
$(\alpha_t,g_t)=(\alpha_1,g_1)$ for $t\leq 1$ and
$(\alpha_t,g_t)=(\alpha_2,g_2)$ for $t\geq 2$. If the homotopy satisfies a
Morse-Smale condition one can now count the number of trajectories of the
time dependent vector field $\pi^*(\grad_{g_t}\alpha_t)$ connecting a critical
point of $\alpha_1$ and a critical point of $\alpha_2$ which have the same
index. This defines a chain mapping
$h_N:C^*_N(\alpha_1,g_1)\to C^*_N(\alpha_2,g_2)$. Let $C_{h_N}^*$ denote
the mapping cone of $h_N$, \cf section~6.

Let $(K,L)$ be a pair of finite CW-complexes, let $\tilde K$ denote the
universal cover of $K$ with the canonic cell structure and let $\tilde L$ be
the preimage of $L$ in $\tilde K$. Then $\pi_1(K)$ acts freely on
$(\tilde K,\tilde L)$ and the cellular cochain complex (with compact
support) $C^*_c(\tilde K,\tilde L)$ \footnote{The component 
$C^k_c(\tilde K,\tilde L)$ consists of the functions with finite support 
defined on the set of $k$-cells of $\tilde K$ which do not belong to 
$\tilde L$.} can be viewed as free finitely
generated $R[\pi_1(K)]$-complex, where $R[\pi_1(K)]$ denotes the group
ring of $\pi_1(K)$ with values in $R$.

We call a Morse-Smale pair $(\alpha,g)$ special if the unstable manifolds
of $\grad_g\alpha$ provide a cell decomposition of $M$.

\proclaim{Proposition 7.1}
Any Morse-Smale homotopy $(\alpha_t,g_t)$ between two special Morse-Smale
pairs $(\alpha_1,g_1)$ and $(\alpha_2,g_2)$ provides a cell decomposition
$(K,L)$ of the compact pair
$\bigl(M\times[-\infty,\infty],M\times\{-\infty\}\bigr)$. Moreover
$$
C^*_{h_N}=\Cal N\otimes_{R[\pi_1(M)]} C^*_c(\tilde K,\tilde L)
$$
as free $\Z$-graded $\Cal N$-complexes. Here $\Cal N$ is regarded as
$R[\pi_1(M)]$-module via the homomorphism $R[\pi_1(M)]\to\Cal N$ induced
by the projection $\pi_1(M)\to\Delta$.
\endproclaim

Now let $(Z,J)$ be a special pair in the sense of Definition~6.3. We then
have an associated special pair $(\alpha,g)$, where $\alpha:=i_Z\omega$ and
$g$ is the Riemannian metric defined by $J$ and $\omega$, \cf section~4.
Let $\iota:M\to\flo$ denote the map that sends points in $M$ to
constant curves. Moreover let $\tilde M\to M$ be a connected component of the
pull back covering $\iota^*\cflo$ and let $\tilde\iota:\tilde M\to\cflo$
be the restriction of the natural map $\iota^*\cflo\to\cflo$. Note that
$\iota_*:\pi_1(M)\to\pi_1(\flo)$ induces a homomorphism
$\Delta\to\Gamma$, where $\Delta$ \resp $\Gamma$
is the structure group of $\tilde M$ \resp $\cflo$, and $\tilde\iota$ is
equivariant.

Since $Z\in\Cal Y^\beta_{\text{$\gamma_0$-reg}}\subset\vfboreg$ every
$1$-periodic trajectory of $Z$ is constant, \ie we have a bijection
$\iota:C\to\Cal P_{\gamma_0}$. Moreover the Conley-Zehnder index
of a constant trajectory equals the Morse index modulo $2$ and hence we
obtain an isomorphism of $\Z_2$-graded $\Lambda$-modules
$$
C^*_F(Z,J)=\Lambda\otimes_{\Cal N}C^*_N(\alpha,g).
\tag{7.1}
$$
As in \cite{FO99}
\footnote {At this point one requires the work of Fukaya-Ono even in the case
of weakly monotone manifolds.} one shows that this actually is an isomorphism
of cochain complexes, although there might be time dependent tubes connecting
two points in $C$, but they do not contribute to the differential of
$C^*_F(Z,J)$.

Suppose $(Z_1,J_1)$ and $(Z_2,J_2)$ be two special pairs and let
$(\alpha_1,g_1)$ and $(\alpha_2,g_2)$ denote the associated special
Morse-Smale pairs. Moreover let $(\alpha_t,g_t)$ be a Morse-Smale homotopy
from $(\alpha_1,g_1)$ to $(\alpha_2,g_2)$ and let
$$
h_N:C^*_N(\alpha_1,g_1)\to C^*_N(\alpha_2,g_2)
\quad\text{\resp}\quad
h_F:C^*_F(Z_1,J_1)\to C^*_F(Z_2,J_2)
$$
denote the chain mappings induced by $(\alpha_t,g_t)$ \resp $(\alpha_t,g_t)$
regarded as homotopy between $(Z_1,J_1)$ to $(Z_2,J_2)$. As in \cite{FO99}
one shows that up to the isomorphism \thetag{7.1} one has
$$
h_F=\Lambda\otimes_{\Cal N}h_N.
$$
Again, there might be time dependent tubes connecting points in $C_1$ and
$C_2$, but they do not contribute to $h_F$. So we have shown

\proclaim{Proposition 7.2}
In the situation above one has for the mapping cone
$$
C^*_{h_F}=\Lambda\otimes_{\Cal N}C^*_{h_N}
$$
as free $\Z_2$-graded $\Lambda$-complexes.
\endproclaim

\demo{Actual proof of Proposition 6.5}
We continue to use the notation introduced above.
From Proposition~7.1 and 7.2 we obtain
$$
C^*_{h_F}=\Lambda\otimes_{R[\pi_1(M)]}C^*_c(\tilde K,\tilde L)
\tag{7.2}
$$
as free $\Z_2$-graded $\Lambda$-complexes, where $(K,L)$ is a CW-complex with
underlying space $\bigl(M\times[-\infty,\infty],M\times\{-\infty\}\bigr)$ and
the $R[\pi_1(M)]$-module structure on $\Lambda$ is given by the composition
$\psi:R[\pi_1(M)]\to\Cal N\to\Lambda$. Since $L$ is a deformation retract of
$K$ one has a Whitehead torsion $\tau(K,L)\in\Wh(\pi_1(M))$, \cf \cite{Mi66}.
From \thetag{7.2} one obtains
$$
\tau\bigl(\gamma_0,(Z_1,J_1),(Z_2,J_2)\bigr)=
\Wh(\psi)\bigl(\tau(K,L)\bigr)\in\Wh(\Lambda),
\tag{7.3}
$$
where $\Wh(\psi):\Wh(\pi_1(M))\to\Wh(\Lambda)$ is the homomorphism induced
from $\psi$. The topological invariance for Whitehead  torsion in the case
compact manifolds
states that $\tau(K,L)$ does only depend on the total space
$\bigl(M\times[-\infty,\infty],M\times\{-\infty\}\bigr)$ and not on the
cell structure. Using a cylindrical cell structure it follows easily,
that $\tau(K,L)=0\in\Wh(\pi_1(M))$. So by \thetag{7.3} the relative torsion
$\tau\bigl(\gamma_0,(Z_1,J_1),(Z_2,J_2)\bigr)$ vanishes too.
\qed
\enddemo

\head
8. Appendix (An example)
\endhead

Consider the torus $T=S^1\times S^1$ with the standard metric $dx^2+dy^2$,
the standard symplectic form $dx\wedge dy$ and the standard complex
structure $J$. Choose bump functions $\lambda$ and $\nu$, such that
$\nu(0)=1$, $\nu'(0)=0$, $\nu''(0)=-1$ and such that $\lambda'$ looks like
\centerline{\epsfbox{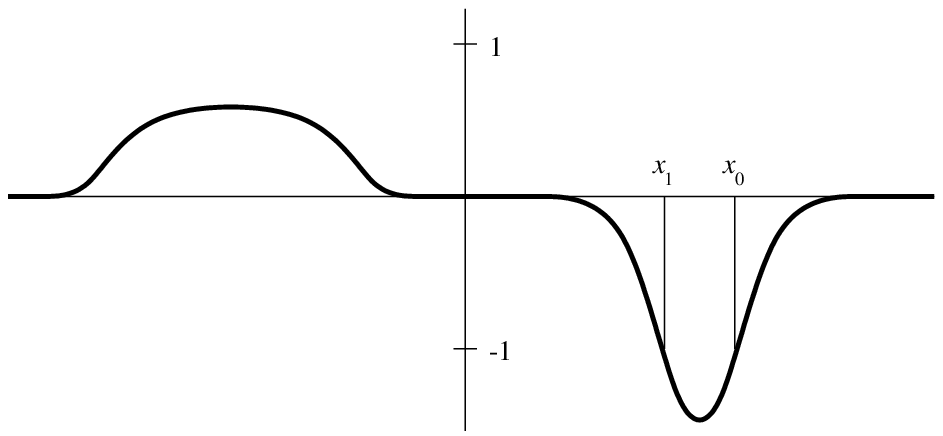}}
where $\lambda'(x_0)=\lambda'(x_1)=-1$. We also assume, that
$0<|\lambda''(x_i)|<2\pi$ and $0<|\lambda(x_i)|<2\pi$.
Consider the time dependent Hamiltonian $h_t(x,y):=\lambda(x)\nu(y-t)$ and
the corresponding time dependent Hamiltonian vector field
$$
Z_t=\lambda(x)\nu'(y-t)\partial_x-\lambda'(x)\nu(y-t)\partial_y.
$$
One easily sees, that the closed curves $t\mapsto(x_i,t)$, $i=0,1$ are
non-degenerate 1-periodic solutions of index difference 1. Indeed, the
index in this example is the Maslov index of the path
$$
t\mapsto e^{t
\Big(\smallmatrix
0                     & -\lambda(x_i)
\\
-\lambda''(x_i) & 0
\endsmallmatrix\Big)
}
$$
which is the Morse index of $J\Big(\smallmatrix 0 & -\lambda(x_i) \\
-\lambda''(x_i) & 0\endsmallmatrix\Big)$, since we assumed
that $|\lambda''(x_i)|<2\pi$, $|\lambda(x_i)|<2\pi$, and this Morse index
is 1 or 2, depending on the sign of $\lambda''(x_i)$.
These should be all periodic solutions in this homotopy class.
However there
are lots of degenerate contractible solutions. Moreover
we have two connecting orbits $u(s,t)=(x(s),t)$,
where $x$ is one of the two (up to shift in $s$) non-constant solutions
of $x'(s)=1+\lambda'(x(s))$. There should not exist other connecting
orbits. Let
$U:=\ker\phi$, $R=\R$. Then $\Lambda_0=\Lambda$ is the field of Laurent series in
one variable, say $z$, and the corresponding $\Z$-graded Floer complex is:
$$
\cdots\to0\to\Lambda@>{1\pm z}>>\Lambda\to0\to\cdots
$$
So its torsion is the non-trivial element
$1\pm z\in\frac{\Lambda\setminus\{0\}}{\pm z}$.

\enddocument